\documentclass[12pt]{article}
\usepackage{amsmath,amsthm,amssymb}
\begin{document}
\title{The suppression thesis in Riemann's Habilitationsvortrag}

\author{Victor Tapia\\ Via dei Giardini 6, 34146 Trieste, Italy\\ tapiens@gmail.com}

\maketitle

\abstract{In 1854 Riemann presented his Habilitationsvortrag {\it \"Uber die Hypothesen welche der Geometrie zu Grunde liegen} to the Faculty of Philosophy of the University of G\"ottingen. The lecture does not contain analytical developments and Riemann never explained the reason for this choice. According to Dedekind this choice was done in order to make the lecture understandable to the non-mathematician professors of the Faculty. This is the suppression thesis. 

Nevertheless, in 1851, three years before, Riemann has presented his Inauguraldissertation {\it Grundlagen f\"ur eine allgemeine Theorie der Functionen einer ver\"anderlichen complexen Gr\"osse} to the same professors of the Faculty of Philosophy. The Inauguraldissertation contains many analytical developments but, on that occasion, they were not suppressed. Therefore, the presence of non-mathematician professors is unrelated with the presence or absence of analytical developments in the lectures.

The ideas contained in Riemann's Habilitationsvortrag --manifold, line element, curvature-- were so new at the time that there was not even a proper language to express them. Therefore, Riemann was forced to adopt a non-technical expository style. This explanation for the absence of analytical developments in Riemann's Habilitationsvortrag is more credible than the suppression thesis.}

\section{Introduction}

On June 10, 1854 Bernhard Riemann presented his Habilitationsvortrag {\it \"Uber die Hypothesen welche der Geometrie zu Grunde liegen} {\it [On the hypotheses which lie at the bases of geometry]} to the Faculty of Philosophy of the University of G\"ottingen \cite{Riemann1854b}. This thesis marks the beginning of modern differential geometry, Riemannian geometry and other developments and applications, such as general relativity. Even when this is a thesis on mathematics, it does not contains analytical developments and Riemann never explained the reason for this choice.

Riemann was born in 1826 and in 1846 he enrolled at the University of G\"ottingen to study philosophy and theology. In 1847 he moved to the University of Berlin to study mathematics and went back to G\"ottingen in 1849 in order to complete his studies of mathematics. Riemann died in 1866 and his {\it Gesammelte Werke} {\it [Collected Works]} \cite{Riemann1876} were published in 1876 by Heinrich Weber as chief editor with the collaboration of Richard Dedekind. A biography of Riemann, the {\it Lebenslauf} {\it [Curriculum Vitae]} \cite{Dedekind1876b}, was prepared by Dedekind and was included in the {\it Gesammelte Werke}.

At the middle of the 19th century the Faculty of Philosophy at the G\"ottingen University comprised not only philosophy but, for administrative reasons, other disciplines. Philology was the most important, but physics and mathematics were also included. According to Dedekind, the non-mathematician professors of the Faculty of Philosophy were not prepared to understand a mathematical exposition. Therefore, Riemann opted for a less demanding expository style with no analytical developments. This hypothetical self-censure on the part of Riemann was named the {\it suppression thesis} by Nowak \cite{Nowak1989}. 

The suppression thesis has become the standard explanation for the absence of analytical developments in Riemann's Habilitationsvortrag. For the benefit of the reader, in an appendix we reproduce extracts of articles and books in which the suppression thesis is mentioned.

In 1851, three years before the Habilitationsvortrag, Riemann has presented his Inauguraldissertation {\it Grundlagen f\"ur eine allgemeine Theorie der Functionen einer ver\"anderlichen complexen Gr\"osse} {\it [Foundations for a general theory of functions of one complex variable]} \cite{Riemann1851} to the same professors of the Faculty of Philosophy. The Inauguraldissertation contains many mathematical developments. On that occasion, however, in spite of the presence of non-mathematician professors, the mathematical developments were not suppressed. 

If the suppression thesis were true, then the Inauguraldissertation would not contain analytical developments. This is not the case. Therefore, the suppression thesis is false. The presence of non-mathematician professors at the Inauguraldissertation and at the  Habilitationsvortrag is unrelated with the presence or absence of analytical developments in those lectures. 

In the Habilitationsvortrag Riemann was introducing new mathematical concepts --manifold, line element, curvature-- which were so new that there was not even a proper language to express them. Therefore, Riemann was forced to adopt a non-technical expository style. This explanation for the absence of analytical developments in Riemann's Habilitationsvortrag is more credible than the suppression thesis.

To support our claim we analyse the concepts of manifold, line element and curvature, as introduced by Riemann, to conclude that they represented departures from previous standard concepts.

\section{The Habilitationsvortrag}

During the 19th century, the academic degrees granted by German universities were the {\it doctorate} and the {\it habilitation}. The doctorate was intended to determine the research abilities of the candidate and was granted after the submission of a {\it Inauguraldissertation} {\it [Doctoral dissertation]}, a written document. The {\it Habilitation} was similar to the {\it venia legendi}, the license to teach. It consisted of two parts: the {\it Habilitationsschrift} {\it [Habilitation thesis]}, a written document submitted to determine the abilities of the candidate for independent research; the second part, the {\it Habilitationsvortrag} {\it [Habilitation lecture]}, was a trial lesson, intended to evaluate the teaching abilities of the candidate. The tradition was that the candidate submitted three titles and the examination committe choose one among them.

The {\it Inauguraldissertation}, the {\it Habilitationsschrift} and the {\it Habilitationsvortrag}, in any discipline, were to be presented to the Faculty of Philosophy, meaning all the professors of the Faculty, even those from other disciplines \cite{Weber1990}. 

In November 1851 Riemann submitted his Inauguraldissertation {\it Grundlagen f\"ur eine allgemeine Theorie der Functionen einer ver\"anderlichen complexen Gr\"osse} \cite{Riemann1851} and the examination took place in December 1851.

The next step was being {\it habilitated} as a {\it Privatdozent} {\it [Private lecturer]}. Accordingly, at the end of 1853 Riemann submitted his Habilitationsschrift {\it \"Uber die Darstellbarkeit einer Function durch eine trigonometrische Reihe} {\it [Concerning the representability of a function by a trigonometric series]} \cite{Riemann1854a}. In December 1853 Riemann submitted three titles for his Habilitationsvortrag: one in analysis, one in algebra and one in geometry. Gauss choose the third topic. Riemann, in a letter to his brother Wilhelm, dated December 28, 1853, wrote \cite{Dedekind1876b}:\footnote{\ All excerpts from the {\it Lebenslauf} are taken from the English translation in \cite{Riemann2004}.}

\begin{itemize}

\item[]{\it My work is now in a reasonable state; I handed in my habilitation thesis at the beginning of December, and I must now propose three subjects for the trial lecture, one of which is then chosen by the faculty. I had already prepared the first two, and I had hoped that one of these would be chosen, but Gauss chose the third, and so I am again in something of a tight spot, as I now have to do some more work on it. I had resumed my other investigations into the connection between electricity, galvanism, light and gravity immediately after finishing my habilitation thesis, and have now made sufficient progress in them, that I can publish my conclusions in this form without hesitation. At the same time I have become more and more convinced that Gauss has been working on this question for several years, and has told a few of his friends -- among others [Wilhelm] Weber -- under the seal of secrecy. Of course I can say this to you, without any risk of appearing presumptuous; anyhow I hope that I am not too late, and that it will be recognised that my results were entirely independent work.}

\end{itemize}

The choice of the third topic by Gauss has been a subject of discussion and speculation. Nevertheless, according to Riemann himself, there was no mystery in this choice and, furthermore, it seems it was not traumatic for Riemann as some authors have claimed. It is interesting to quote some paragraphs of Laugwitz \cite{Laugwitz1999} since they illustrate the lack of agreement existing with respect to the choice of the third topic. He writes:

\begin{itemize}

\item[]{\it As is done even today, Riemann was supposed to submit titles of three lectures, from which the Faculty of Philosophy was to choose one. The first title was ``Geschichte der Frage \"uber die Darstellbarkeit einer Function durch eine trigonometrische Reihe'' [The history of the question of representability of a function as a trigonometric series], and the second, ``\"Uber die Aufl\"osung zweier Gleichungen zweiten Grades mit zwei unbekannten Gr\"ossen'' [Solution of two quadratic equations in two unknowns]. Of course, the third title was ``\"Uber die Hypothesen $\cdots$''. To this day it is expected that the candidate will propose topics that belong to very different areas and thus demonstrate the breadth of his professional knowledge. Accordingly, Riemann proposed topics in analysis, algebra, and geometry. The Faculty's, in this case Gauss, choice was not, as Dedekind put it, a breach ``of the usual convention'' and his suggestion that Gauss ``chose the third [topic] because he was curious to hear how such a young person would handle so difficult a topic'' was pure guesswork. The first topic was largely taken care of by the dissertation [the  Habilitationsschrift], and the second, the determination of the points of intersection of two quadratic curves, must have struck Gauss as trivial. Even if we take seriously the suggestion that someone of Riemann's stature could change even unpromising material into something unusual, the Faculty could not but reject the first two topics. Riemann could not have been naive enough not to know this, and there is every reason to think that he planned to lecture about the ideas that were particularly important to him.}

\item[]{\it There are statements in the literature that Riemann expressed great surprise at the Faculty's choice of the third topic. He did not. By December 1853 he had probably worked out the details of the first two topics but not of the third. No manuscripts bearing on the elaboration of any of the three have been preserved. Riemann had been assigned the third topic as early as December. He seems not to have been overly worried by it, for he postponed its elaboration until Easter of 1854. This seems to indicate that he knew very well what he wanted to say in his lecture. At that time he was always occupied with issues of physical and philosophical relations. Mathematical issues, the concept of a manifold and Gauss' geometry of surfaces, had been touched upon in his dissertation. It is impossible to tell when the idea of an $n$-dimensional generalization occurred to him, and whether this was the result of his preoccupation with Herbart. Completion of the computational part of his work was to cause him fewer difficulties than preparation of a polished presentation.}

\end{itemize}

The Habilitationsvortrag took place on June 10, 1854. The best description is given by Riemann himself in a letter to his brother, dated June 26, 1854 \cite{Dedekind1876b,Gallagher1981}:

\begin{itemize}

\item[]{\it Around Christmas, if I remember rightly, I wrote to you from G\"ottingen, telling you that I had completed my habilitation thesis at the beginning of December, and handed it in to the Dean, and that immediately afterwards I had continued with my researches into the interrelations between the fundamental laws of physics and had become so engrossed in this study that when the subject was set for the trial lecture at the colloquium, I could not tear myself away from it. Soon after this I became ill, partly as a result of too much brooding, and partly because of sitting in front of a stove during this foul weather; my old trouble started again with a vengeance, which was not without adverse effect on my work. It was not until several weeks later, when the weather improved and I was able to get about a bit, that my health got better. I have now rented a summer-house for the summer and since then, thank god, have had nothing to complain of as regards my health. For about a fortnight after Easter [Easter Sunday 16th April] I had to deal with other work which I could not very well avoid. That done, I went ahead enthusiastically with the preparation of my trial lecture and had completed it by Pentecost [Sunday, June 4]. Although it was quite an effort, I needed to have it done so that I could have my examination immediately, and would not have to visit Quickborn without having achieved my object. In fact, Gauss's state of health has recently become so bad, that everyone now fears that he will die this year,\footnote{\ Gauss passed away on February 23, 1855.} and he himself felt that he was too weak to examine me. He now asked me to wait since in any case I could not start lecturing until the next semester -- at least until August, when he might be feeling better. I had already resolved to bow to the inevitable. Then suddenly at midday on the Friday after Pentecost [June 9] he decided to give in to my repeated pleas to get rid of ``the halter round my neck'' by agreeing to hold the examination at half past ten on the following day [Saturday, June 10], and so by one o'clock on the Saturday it was all happily ended.}

\end{itemize}

\noindent Riemann does not give any further information on the circumstances of the habilitation lecture neither any clues on why he choose a non-technical expository style.

Riemann died in 1866 and his wife, Elise Koch, asked Dedekind to take care of the publication of Riemann's mathematical legacy \cite{Dedekind1876a}. In 1867 Dedekind published three of Riemann's manuscripts: {\it Ein Beitrag zur Elektrodynamik} {\it [A contribution to electrodynamics]} \cite{Riemann1867}; the Habilitationsschrift {\it \"Uber die Darstellbarkeit einer Function durch eine trigonometrische Reihe} {\it [Concerning the representability of a function by a trigonometric series]} \cite{Riemann1868a}; and the Habilitationsvortrag {\it \"Uber die Hypothesen, welche der Geometrie zu Grunde liegen} \cite{Riemann1868b}.

For the rest of Riemann's mathematical production, in 1872 Dedekind looked for the help of Alfred Clebsch and it was agreed that Clebsch would act as the main editor, but he died in November of 1872. In 1874 Dedekind contacted H. Weber to undertake the rest of the edition. Riemann's {\it Gesammelte Werke} were published in 1876 \cite{Riemann1876} with H. Weber as chief editor and Dedekind appearing as a collaborator. A biography of Riemann, the {\it Lebenslauf} \cite{Dedekind1876b}, was prepared by Dedekind, the {\it Lebenslauf}, and included in Riemann's {\it Gesammelte Werke}.

\section{The suppression thesis}

According to Dedekind \cite{Dedekind1876b}, Riemann wanted to make the Habilitationsvortrag understandable to the non-mathematician professors of the Faculty of Philosophy, and this would explain the absence of analytical developments:

\begin{itemize}

\item[]{\it As regards the first part of this letter, it should be noted that Riemann, in composing his trial lecture on the hypotheses of geometry, had made his task very much harder by striving to ensure that it would be as easily understandable as possible even to members of the Faculty who had no mathematical training. By doing so, however, he created a most admirable masterpiece of exposition, inasmuch as, while suppressing all the detailed mathematical analysis, he nevertheless succeeds in conveying so accurately the train of his thoughts, that it can be completely reconstructed from the indications which he gave. Gauss, contrary to the usual practice, had selected from the three proposed themes not the first but the third, because he was keen to know how such a difficult subject would be handled by such a young man. The lecture, which exceeded all his expectations, greatly astonished him, and on his way back from the meeting of the faculty he spoke to Wilhem Weber enthusiastically, and quite uncharacteristic excitement, about the profundity of the ideas which Riemann had put forward.}

\end{itemize}

\noindent This hypothetical self-censure on the part of Riemann is the {\it suppression thesis}. The term was first used in this respect by Nowak \cite{Nowak1989}.

The suppression thesis is a statement found in the {\it Lebenslauf} prepared by Dedekind. That statement was not made by Riemann. Among the scarse surviving documentation there is no mention by Riemann himself of such decision.

In December 1851, three years before the Habilitationsvortrag, Riemann presented his Inauguraldissertation {\it Grundlagen f\"ur eine allgemeine Theorie der Functionen einer ver\"anderlichen complexen Gr\"osse} \cite{Riemann1851} to the same professors of the Faculty of Philosophy. The Inauguraldissertation contains many analytical developments but, on that occasion, they were not suppressed. Therefore, the presence of non-mathematician professors at the lectures is unrelated to the presence or absence of analytical developments in the lectures. 

If the suppression thesis were true, then the Inauguraldissertation would not contain analytical developments. This is not the case. Therefore, the suppression thesis is false.

\section{The content of the Habilitationsvortrag}

The new ideas introduced by Riemann in his Habilitationsvortrag were: {\it manifold}, {\it line element} and {\it curvature}. They were so new that there was not even a proper language to express them. An algorithmic language was inexistent. Therefore, Riemann was forced to adopt a non-technical expository style.

In 1828 Gauss published his {\it Disquisitiones generales circa superficies curvas} {\it [General considerations on curved surfaces]} \cite{Gauss1828}. In this work surfaces are considered as two-dimensional extensions embedded in a three-dimensional space. Gauss introduces the curvature of a surface as a measure of how much the curved surface departs from flatness.

The other important result in the {\it Disquisitiones} is the {\it Theorema Egregium}, the fact that the curvature of a surface depends only on intrinsic properties, that is, it can be determined by means of measurements of distances on the surface. In this way, there is no reference to the way in which the surface is embedded in a three-dimensional space. 

Since the geometry of a surface does not depend anymore on the ambient three-dimensional space, less on its dimension, there are no restrictions on the dimensionality of the surface. We are now free to consider `surfaces' of any dimension. That was the step taken by Riemann in his habilitation lecture. The concept of surface, as a two-dimensional extension embedded in three dimensions, was replaced by that of an $n$ dimensional extension with no reference to an embedding space.

In Euclidean geometry a {\it line} is a one-dimensional extension that lies evenly with respect to all its points and can be prolonged indefinitely in each direction; a {\it plane} is a two-dimensional extension that lies evenly with respect to all its points and can be prolonged indefinitely in all directions. If we lift the restriction of `evenly' we arrive to the concepts of {\it curve} and {\it surface}: a {\it curve} is a one-dimensional extension which can be prolonged in both directions; a {\it surface} is a two-dimensional extension which can be prolonged in all directions.

A {\it manifold} of dimension $n$ is an extension which can be prolonged along all of its $n$ directions. A manifold can be prolonged in many different ways, and this gives rise to manifolds with different topology. {\it Differentiable manifolds} are manifolds with a smooth structure and then it is possible to introduce coordinates.

In order to express his idea of an $n$ dimensional extension, Riemann coined the term ``Mannigfaltigkeit''; Clifford \cite{Clifford1873} translated it as ``manifoldness'', a word inexistent in English, a clear evidence that this was a new word. Later on it became ``manifold''.

The {\it line element} is a measure of the infinitesimal distance betwen two points along a given curve. The line element $ds$ depends on the position, that is, on the coordinates, and the direction in which the distance is measured, that is, on the differentials of the coordinates. If the curve is divided into infinitesimal pieces and each infinitesimal piece is increased in the same proportion, then the distance must increase proportionally. From here Riemann concludes that the line element is a function which is homogeneous of the first order in the differentials of the coordinates.

The simplest solution occurs when the line element is the square root of an expression quadratic in the differentials of the coordinates. The concept of Euclidean distance, based on the Pythagorean Theorem, was replaced by that of a line element quadratic in the differentials of the coordinates. This is Riemannian geometry and the coefficients in the quadratic differential form are the components of the {\it metric} tensor.

The {\it curvature} of higher-dimensional spaces is obtained as a generalization of the curvature of two-dimensional surfaces. In the $n$-dimensional space we consider the $n\,(n-1)/2$ possible two-dimensional subspaces. For each of these two-dimensional subspaces we construct the Gaussian curvature. Then, the curvature of a higher-dimensional space is the set of all possible Gaussian (two-dimensional) curvatures.

During the second half of 19th century tensor calculus was developed. The final form of tensor calculus was developed by Ricci and Levi-Civita in 1900. With tensor calculus as a tool it was then possible to give a reasonably modern presentation of Riemannian geometry. The sectional curvatures introduced by Riemann are the components of the Riemann-Christoffel tensor.

In this way we see that the three main ideas considered by Riemann --manifold, line element, curvature-- departed considerably from the previously corresponding concepts.

\section{Summary and conclusions}

Riemann's Habilitationsvortrag {\it \"Uber die Hypothesen welche der Geometrie zu Grunde liegen} does not contain analytical developments. According to Dedekind, Riemann deliberately choose a non-technical expository style in order to make his thesis understandable to the non-mathematician professors of the Faculty of Philosophy. This is the suppression thesis. However, Riemann's Inauguraldissertation, presented three years before to the same professors of the Faculty of Philosophy, contains mathematical developments which were not suppressed. The conclusion is that the suppression thesis is false, it is an untenable explanation for the absence of analytical developments in Riemann's Habilitationsvortrag. 

The ideas presented in Riemann's Habilitationsvortrag --manifold, line element, curvature-- were so new that there was not even a proper language to express them. Therefore, Riemann was forced to adopt a non-technical expository style. This is a more credible explanation for the absence of analytical developments in the Habilitationsvortrag than the suppression thesis.

To conclude, the suppression thesis is a myth established just by force of repetition.

\section*{Appendix. The suppression thesis in the literature}

The suppression thesis has become the standard explanation for the absence of analytical developments in the Habiltationsvortrag and has been mentioned in several works. For the benefit of the reader we reproduce extracts of articles and books in which the suppression thesis is mentioned.

In order to help trace the transmission of the suppression thesis, at the end of each extract we cite the works in which it is based.\footnote{\ Boldface words or phrases are my emphasis and refer to statements of dubious validity. Square brackets are used for extrapolations.}

\begin{itemize}

\item[$\bullet$]Bonola (1906, 1912) \cite{Bonola1906}:

\begin{itemize}

\item[]{\it Fu letta da Riemann nel 1854, per la sua abilitazione presso la Facolt\`a filosofica di Gottinga, davanti ad un pubblico composto non di soli matematici. Perci\`o non contiene sviluppi analitici ed i concetti ivi esposti hanno veste prevalentemente intuitiva.}

\item[]{\it It was read by Riemann to the Philosophical Faculty at G\"ottingen as his Habilitationsschrift, before an audience not composed solely of mathematicans. For this reason it does not contain analytical developments, and the conceptions introduced are mostly of an intuitive character.}

\item[$\bullet$]Dedekind is cited, not in connection with the {\it Lebenslauf} but as an editor of Riemann's {\it Gesammelte Werke} \cite{Dedekind1876b}.

\end{itemize}

\item[$\bullet$]Spivak (1970) \cite{Spivak1970}:

\begin{itemize}

\item[]{\it Riemann hoped to make his lecture intelligible even to those members of the Faculty who knew little mathematics. Consequently, hardly any formulas appear and the analytic investigations are completely suppressed. Although Dedekind describes the lecture as a masterpiece of exposition, it is questionable how many of the Faculty comprehended it.}

\item[$\bullet$]Dedekind (1876) \cite{Dedekind1876b}.

\end{itemize}

\item[$\bullet$]Kline (1972) \cite{Kline1972}.

\begin{itemize}

\item[]{\it Before examining the details we should be forewarned that Riemann's ideas as expressed in the lecture and in the manuscript of 1854 are vague. One reason is that Riemann adapted it to his audience, the entire Faculty at G\"ottingen. Part of the vagueness stems from the philosophical considerations with which Riemann began his paper.}

\item[$\bullet$]No reference for the suppression thesis is given.

\end{itemize}

\item[$\bullet$]Portnoy (1982) \cite{Portnoy1982}.

\begin{itemize}

\item[]{\it $\cdots$ Riemann had worked hard to make the lecture understandable to nonmathematicians in the audience, $\cdots$}

\item[]{\it Riemann's definition [of manifold] is vague and awkward by comparison [to modern ones], but it has an important advantage in being constructive rather than analytic $\cdots$. We may regard his choice of an intuitive definition as an indication of the desire to be intelligible to his general scholarly audience.}

\item[$\bullet$]Dedekind (1876) and Spivak (1970) are cited but no direct reference is found in the text.

\end{itemize}

\item[$\bullet$]Reich (1994) \cite{Reich1994}.

\begin{itemize}

\item[]{\it Riemann h\"alt seinen Habilitationsvortrag nicht in erster Linie f\"ur Mathematiker, er ist vielmehr bem\"uht, seine Gedanken mit so wenig Mathematik wie m\"oglich zu pr\"asentieren.}

\item[]{\it Riemann does not give his habilitation lecture primarily for mathematicians; rather, he tries to convey his thoughts with as little mathematics as possible to present.}

\item[$\bullet$]Nowak (1989) \cite{Nowak1989}.

\end{itemize}

\item[$\bullet$]Gray (2007) \cite{Gray2007}.

\begin{itemize}

\item[]{\it It lacks the formulae which would help mathematicians understand it, because the lecture was given to the Philosophy Faculty, of which mathematics formed a department.}

\item[$\bullet$]Dedekind (1876) \cite{Dedekind1876b}.

\end{itemize}

\item[$\bullet$]Pesic (2007) \cite{Pesic2007}.

\begin{itemize}

\item[]{\it For Riemann {\bf con­sciously} chose to express his ideas with only a single equation, ``almost without calculation,'' {\bf as he put it}, the better to bring forward their philosophical import.}

\item[$\bullet$]Spivak (1970) \cite{Spivak1970}, Portnoy (1982) \cite{Portnoy1982}, Nowak (1989) \cite{Nowak1989}.

\end{itemize}

\item[$\bullet$]Giovanelli (2013) \cite{Giovanelli2013}.

\begin{itemize}

\item[]{\it Riemann's lecture, which was intended for an audience of non-mathema\-ticians, was {\bf intentionally} scarce in the use of mathematical formulas.}

\item[$\bullet$]The {\it Lebenslauf} is not cited.

\end{itemize}

\item[$\bullet$]Darrigol (2015) \cite{Darrigol2015}.

\begin{itemize}

\item[]{\it Riemann was addressing a mostly non-mathematical\ \ audience,\ \ the G\"ottingen Faculty of Philosophy, who {\bf would not have} followed the mathematical technicalities.}

\item[$\bullet$]Dedekind (1876) \cite{Dedekind1876b}, Portnoy (1982) \cite{Portnoy1982}, Gray (2007) \cite{Gray2007}.

\end{itemize}

\item[$\bullet$]Jost (2013) \cite{Jost2013}.

\begin{itemize}

\item[]{\it It is noteworthy that Riemann's ``Hypotheses'' as one of the key texts in mathematics proceeds without mathematical formulas (in the whole text, there is only a single formula which is of only marginal importance). This sets Riemann's text apart from other foundational mathematical works, like the sophisticated and deeply thought out symbolism of Leibniz or the formalization of the infinite of Cantor. Even his most important precursor, Gauss' ``Disquisitiones generales circa superficies curvas'', which founded modern differential geometry, the starting point of Riemannian Geometry, is different in this respect. At least in this case, the history of mathematics is not simply a progressive formalization, but it turns out that mathematical abstraction can in principle rise well above formulas.}

\item[]{\it Of course, the occasion for Riemann's paper must also be considered: it was a colloquium before the Faculty of Philosophy, and Riemann {\bf certainly} wanted to take into account the lack of mathematical expertise of most of his listeners. Among them, besides Gauss, who, incidentally, was not a professor of mathematics but of astronomy and director of the Observatory, mathematics was only represented by the two professors Ulrich (1798-1879) and Stern (1807-1894). However, other similar lectures and writings, such as Klein's Erlangen Program, with which he [Klein] introduced himself to the faculty in Erlangen, certainly display a much more formulaic character, and if the Faculty [G\"ottingen] had chosen one of the other topics proposed by Riemann, the presentation would presumably have been developed in mathematical formulas as well.}

\item[$\bullet$]Dedekind (1876) \cite{Dedekind1876b}, Spivak (1970) \cite{Spivak1970}.

\end{itemize}

\item[$\bullet$]Gray (2017) \cite{Gray2017}.

\begin{itemize}

\item[]{\it $\cdots$ it was given to the Philosophy Faculty at G\"ottingen, of which Mathematics was a Department, with Gauss as one of the examiners. These circumstances explain the unfortunate absence of formulae that would otherwise have assisted subsequent readers.}

\item[$\bullet$]Dedekind (1876) \cite{Dedekind1876b}, Spivak (1970) \cite{Spivak1970}, Darrigol (2015) \cite{Darrigol2015}.

\end{itemize}

\item[$\bullet$]Ji (2017) \cite{Ji2017}.

\begin{itemize}

\item[]{\it Since it was the written version of an address towards the general Faculty of G\"ottingen University, it is not a technical mathematical paper, and does not contain formulas and precise definitions.}

\item[]{\it Riemann {\bf tried} to prepare his lecture so that members of the Faculty without mathematical training could understand it and hence suppressed all the detailed mathematical computation.}

\item[$\bullet$]No reference for the suppression thesis is given.

\end{itemize}

\end{itemize}

Alternative explanations for the absence of analytical developments in Riemann's Habilitationsvortrag have also been offered. Among them we have:

\begin{itemize}

\item[$\bullet$]Kagan (1949) \cite{Kagan1949}:

\begin{itemize}

\item[]{\it According to Dedekind, the lecture was written so as to be accessible to all professors of the Faculty, most of whom knew no mathematics. This is hardly credible. What is far more likely is that Riemann prepared a lecture for Gauss. This explains the form of the presentation. In a few pages Riemann presents, or rather sketches, a number of profound ideas stated in strictly mathematical terms. No computations are given -- only results, presented in a vague and compressed form. This is a piece of unfinished research rather than a memoir prepared by an author for publication. Riemann regarded his lecture as insufficiently elaborated and did not publish it. After Riemann's death, Dedekind extracted the manuscript from his ``Nachlass'' and published it in 1868.}

\item[$\bullet$]Dedekind is cited in the text, but not in the references.

\end{itemize}

\item[$\bullet$]Lewy (1953) \cite{Lewy1953}.

\begin{itemize}

\item[]{\it Avoiding the help of formulae in order to stress the non-technical character of his ideas, $\cdots$.}

\end{itemize}

\item[$\bullet$]Nowak (1989) \cite{Nowak1989}:

\begin{itemize}

\item[]{\it A look at {\bf Riemann's audience}, and those he cites as influences, is also relevant to a discussion of his intent in writing the ``Habilitationsvortrag''. It is worth noting that Gauss chose the third of a list of three topics because he was interested in what Riemann would have to say on such a difficult topic; Riemann had prepared the first two topics but not the third. Although Riemann had recorded some thoughts about relations of extension earlier in his life, so far as we know he never organized this material in any form. We may conclude that except for the mathematical ideas used, Riemann wrote the paper in a few weeks during the spring of 1854, fully aware that it was to be delivered as a lecture rather than printed. He referred to it beforehand as his ``Probevorlesung'' or ``trial lecture'', so we can assume that he had his expected audience -- the Philosophical Faculty of G\"ottingen -- in mind as he wrote. All commentators on the paper, from Dedekind on down, have suggested that Riemann simplified his discussion for the sake of his audience. If we do not interpret this consideration as merely a negative one -- one that prompted him to reduce the mathematical content -- but as a positive consideration -- one that led him to speak to the philosophical interests of a {\bf large part of his audience} -- we can provide a better explanation for the structure and themes of the ``Habilitationsvortrag''. Riemann was demonstrating that he could perform the functions of the ideal Dozent, conducting detailed, original research and relating it to the larger concerns of the culture.}

\item[]{\it The impact Riemann's ``Habilitationsvortrag'' eventually had upon mathematics is unquestionable. But this should not lead one to assume too quickly that Riemann was only interested in presenting mathematics. Some authors have presented a ``suppression'' thesis about the ``Habilitationsvortrag'': that Riemann wished to present detailed mathematics, but could not because there were non-mathematicians in his audience, and so he adopted a vague philosophical standpoint which actively detracted from the clarity of the presentation. Several points, for example Riemann's failure to publish a paper surveying his results in differential geometry, and the rhetorical structure of the ``Habilitationsvortrag'', do not so much support the ``suppression'' thesis as the possibility that Riemann had something other than geometry in mind as the purpose of the ``Habilitationsvortrag''.}

\item[]{\it The suppression thesis also does a disservice to Riemann's reputation as a mathematical expositor: we are asked to believe that Riemann wished to present mathematics, but was forced to render his presentation vague and enigmatic in order to make this mathematics intelligible to nonmathematicians. Even if simplification was required for the sake of the audience, it does not account for either the extended investigation of the concept of space in the first section, or the discussion of applications to space in the third. If, on the contrary, we recall that the ``Habilitationsvortrag'' was the text of a lecture, delivered by a man well read in philosophy to an audience more learned in philosophy than in mathematics, and was part of an occasion whose purpose it was to demonstrate lecturing ability in a culture whose ideal researcher could relate his specific investigations to larger intellectual issues, we could well conclude that the paper was primarily philosophical rather than mathematical in intent. Such a view explains many features of the history of the ``Habilitationsvortrag'' -- the fact that it was never published as a mathematical paper, its unusual organization, which placed most of the mathematical content in one section, and its reception, which took issue with it primarily on philosophical rather than mathematical grounds. In the ``Habilitationsvortrag'', Riemann is functioning as more than a mathematician; we have a better understanding of the paper if we see him speaking primarily as a philosopher, possessed of a rather powerful mathematical methodology. Riemann was tactfully suggesting that new mathematical researches might have something to say to philosophy. His attempt to criticize the Kantian view of space on mathematical grounds presented mathematics as a significant and contributing component of intellectual culture.}

\item[]{\it $\cdots$ we have a better understanding of the paper if we see him [Riemann] speaking primarily as a philosopher, possessed of a rather powerful mathematical methodology.}

\item[$\bullet$]Dedekind (1876) \cite{Dedekind1876b}, Spivak (1970) \cite{Spivak1970}, Portnoy (1982) \cite{Portnoy1982}.

\item[]The term suppression thesis was first used in this article.

\end{itemize}

\item[$\bullet$]Laugwitz (1999) \cite{Laugwitz1999}.

\begin{itemize}

\item[]{\it In the lecture, Riemann pushed to extremes {\bf his tendency} to use as few formulas as possible.}

\item[]{\it Now we see that the absence of formulas in the lecture was not a gesture aiming to accommodate the nonmathematical members of the Faculty. Rather, it was an attempt on Riemann's part to exhibit mathematics, in Herbart's sense, as part of philosophy, as thinking in concepts.}

\item[$\bullet$]No reference is given for the suppression thesis.

\item[]The emphasized words do not have any support. Most of Riemann's works contain analytical developments and some of them quite generously; see Edwards (2010) \cite{Edwards2010}.

\end{itemize}

\item[$\bullet$]Scholz (1999) \cite{Scholz1999}.

\begin{itemize}

\item[]{\it He was completely aware that he was working in a border region between mathematics, physics, and philosophy, not only in the sense of the pragmatic reason that his audience was mixed, but by the very nature of his exposition. There was no linguistic or symbolical frame inside mathematics, which he could refer to, even only to formulate a general concept of manifold}.

\end{itemize}

\item[$\bullet$]Gray (2007) \cite{Gray2007}.

\begin{itemize}

\item[]{\it $\cdots$ and the fact that mathematical terminology was being created and complexities uncovered as Riemann {\bf and Clifford} were writing.}

\item[$\bullet$]Dedekind (1876) \cite{Dedekind1876b}.

\item[]Clifford acted only as translator of Riemann's Habilitationsvortrag in 1873.

\end{itemize}

\end{itemize}


\end{document}